\documentclass[a4paper,english,12pt]{article}
\usepackage[logonly]{trace}
\usepackage{babel}
\usepackage{amsfonts, amsmath, amssymb}
\usepackage{graphicx}
\usepackage{pstricks}
\usepackage{psfig}
\usepackage{multirow}
\usepackage{fancyhdr}
\usepackage{vmargin, fancybox}

\newcommand{\mathsym}[1]{{}}

\usepackage{theorem} \theorembodyfont{\upshape}

\newtheorem{theorem}{Theorem}[section]

\newtheorem{proposition}[theorem]{Proposition}

\newtheorem{definition}[theorem]{Definition}

\usepackage{graphicx}
\usepackage{pstricks}
\usepackage{psfig}

\title{DRP scheme optimization}

\author{Claire David *$\dag$ and Pierre Sagaut \thanks{Universit\'e Pierre et Marie Curie-Paris 6,
Laboratoire de Mod\'elisation en M\'ecanique, UMR CNRS 7607, Bo\^ite
courrier $n^0162$, 4 place Jussieu, 75252 Paris cedex 05, France -
tel. (+33)1.44.27.62.13; Fax (+33)1.44.27.52.59   ($\dag$
corresponding author: { david@lmm.jussieu.fr}).}}

\begin{document}

\maketitle

\begin{abstract}
A new DRP scheme is built, which enables us to minimize the error
due to the finite difference approximation, by means of an
equivalent matrix equation.
\end{abstract}

\noindent{\textbf{keywords}}\\DRP schemes, Sylvester equation

\pagestyle{myheadings} \thispagestyle{plain} \markboth{CL. DAVID AND
P. SAGAUT}{DRP SCHEME OPTIMIZATION}

\section{Introduction: Scheme classes}
\label{sec:intro}

\indent We hereafter propose a method that enables us to build a DRP
scheme while minimizing the error due to the finite difference
approximation, by means of an equivalent matrix equation.\\
\\

\noindent Consider the transport equation:
\begin{equation}
\label{transp} \frac{\partial u}{\partial t}+\frac{\partial
u}{\partial x}=0 \,\,\, , \,\,\, x \,\in \,[0,L], \,\,\,t \,\in
\,[0,T]
\end{equation}

\noindent with the initial condition $u(x,t=0)=u_0(x)$.

\bigskip

\begin{proposition}

\noindent A finite difference scheme for this equation can be
written under the form:
\begin{equation} \label{scheme} {{{{{\alpha \, u}}_i}}^{n+1}}+
   {{{{{\beta \,u}}_i}}^{n}}
    +{{{{{\gamma \,u}}_i}}^{n-1}}
      +\delta \,{{{u_{i+1}}}^n}+{{{{{\varepsilon \, u}}_{i-1}}}^n}
            +{{{{{\zeta \,u}}_{i+1}}}^{n+1}}
              +{{{{{\eta \,u}}_{i-1}}}^{n-1}}+{{{{{\theta \,u}}_{i-1}}}^{n+1}}+\vartheta \,{{ u}_{i+1}}^{n-1} =0
              \end{equation}

\noindent where:
\begin{equation}
{u_l}^m=u\,(l\,h, m\,\tau)
\end{equation}
\noindent  $l\, \in \, \{i-1,\, i, \, i+1\}$, $m \, \in \, \{n-1,\,
n, \, n+1\}$, $j=0, \, ..., \, n_x$, $n=0, \, ..., \, n_t$, $h$,
$\tau$ denoting respectively the mesh size and time step ($L=n_x\,h$, $T=n_t\,\tau$).\\
The Courant-Friedrichs-Lewy number ($cfl$) is defined as $\sigma = c \,\tau / h$ .\\
\\

A numerical scheme is  specified by selecting appropriate values of
the coefficients  $\alpha$, $\beta$, $\gamma$, $\delta$,
$\varepsilon$, $\zeta$, $\eta$, $\theta$ and  $\vartheta$ in
equation (\ref{scheme}), which, for sake of usefulness, will be
written as:
\begin{equation}
\alpha=\alpha_x+\alpha_t\,\,\, , \,\,\,\beta=\beta_x+\beta_t\,\,\, ,
\,\,\,\gamma=\gamma_x+\gamma_t\,\,\, ,
\,\,\,\delta=\delta_x+\delta_t\,\,\, ,
\,\,\,\varepsilon=\varepsilon_x+\varepsilon_t\,\,\, , \,\,\,
\end{equation}
\noindent where the "$_x$" denotes a dependance towards the mesh
size $h$, while the "$_t$" denotes a dependance towards the time
step
$\tau$.\\
\noindent  Values corresponding to numerical schemes retained for
the present works are given in Table \ref{SchemeTable}.
\end{proposition}
\bigskip

\begin{table}
\caption{Numerical scheme coefficient.}
\begin{center} \footnotesize
{\begin{tabular}{cccccccccc} \hline
Name & $  \alpha $ & $ \beta$ & $\gamma$ & $\delta $  & $\epsilon$ & $\zeta$ & $\eta$ & $\theta$ & $\vartheta$ \\
\hline
 & $ \alpha_x+\alpha_t $ & $ \beta_x+\beta_t$ & $\gamma_x+\gamma_t$ & $\delta_x+\delta_t $  &
 $\varepsilon_x+\varepsilon_t$ &  &  &  &  \\
\hline
Leapfrog & $\frac{1}{2 \tau} $ & 0 &  $\frac{-1}{2 \tau} $ & $
\frac{1}{ 2 h}  $ &
$   \frac{-1}{ 2 h}  $ & 0 & 0 & 0  & 0 \\
Lax & $\frac{1}{ \tau} $ &  0 & 0 & $  \frac{1}{ 2 h } - \frac{1}{ 2
\tau} $ &
$  - \frac{1}{ 2 h}- \frac{1}{2 \tau}  $ & 0 & 0 & 0  & 0 \\
Lax-Wendroff & $\frac{1}{ \tau} $ & $  \frac{1 \tau}{ h ^2}-\frac{1}{ \tau}   $ & 0 & $ \frac{1- \sigma }{ 2 h} $
&  $  \frac{-( 1+ \sigma ) }{ 2 h } $ & 0 & 0 & 0 & 0 \\
Crank-Nicolson & $ \frac{1}{ h ^2}+\frac{1}{ \tau}   $ &  $
\frac{1}{ h ^2}-\frac{1}{ \tau}   $ & 0 & $ \frac{-1}{ h ^2} $ & $
\frac{-1}{ h ^2 } $ & 0 & $ \frac{-1}{ h ^2} $ & $ \frac{-1}{ h ^2}
$ & 0
\end{tabular}}
\end{center}
\label{SchemeTable}
\end{table}

\noindent The number of time steps will be denoted $n_t$, the number
of space
steps, $n_x$. In general, $n_x\gg n_t$.\\

\noindent In the following: the only dependance of the coefficients
towards the time step $\tau$ existing only in the Crank-Nicolson
scheme, we will restrain our study to the specific case:

\begin{equation}
 \alpha_t=\gamma_t=\zeta=\eta=\theta=\vartheta=0
\end{equation}

\noindent The paper is organized as follows. The building of the DRP
scheme is exposed in section \ref{DRP}. The equivalent matrix
equation, which enables us to minimize the error due to the finite
difference approximation, is presented in section \ref{Sylv}. A
numerical example is given in section \ref{Ex}.

\section{The DRP scheme}
\label{DRP}

 \noindent The first derivative $\frac{\partial
u}{\partial x}$ is approximated at the $l^{th}$ node of the spatial
mesh by:

\begin{equation}\label{approx}
 (\, \frac{\partial u}{\partial
x}\,)_l  \simeq
   {{{{{\beta_x \,u}}_{l+i}}}^{n}}
      +\delta_x \,{{{u_{l+i+1}}}^n}+{{{{{\varepsilon_x\,
      u}}_{l+i-1}}}^n}
\end{equation}
\noindent Following the method exposed by C. Tam and J. Webb in
\cite{Tam}, the coefficients $\beta_x$, $\delta_x$, and
$\varepsilon_x$ are determined requiring the Fourier Transform of
the finite difference scheme (\ref{approx}) to be a close
approximation of the partial derivative $ (\, \frac{\partial
u}{\partial x}\,
)_l$.\\
\noindent (\ref{approx}) is a special case of:

\begin{equation}\label{approx_Cont}
 (\, \frac{\partial u}{\partial
x}\,)_l  \simeq
   \beta_x \,u(x+i\,h)
      +\delta_x \,u(x+(i+1)\,h)+\varepsilon_x\,u(x+(i-1)\,h)
\end{equation}

\noindent where $x$ is a continuous variable, and can be recovered
setting $x=l\,h$.\\
\noindent Applying the Fourier transform, referred to by
$\,\widehat{\, }$ , to both sides of (\ref{approx_Cont}), yields:

\begin{equation}
\label{Wavenb}
 j\, \omega \, \widehat{u}  \simeq \left \lbrace
   \beta_x \,e^{\,0}
      +\delta_x \,e^{\,j\,\omega\,h}+\varepsilon_x\,e^{\,-\,j\,\omega\,h}
     \right \rbrace \, \widehat{u}
\end{equation}
\noindent  $j$ denoting the complex square root of $-1$.\\




\noindent Comparing the two sides of (\ref{Wavenb}) enables us to
identify the wavenumber $ \overline{\lambda}$ of the finite
difference scheme (\ref{approx}) and the quantity $\frac{1}{j}\,
\left \lbrace\beta_x \,e^{\,0}
      +\delta_x
      \,e^{\,j\,\omega\,h}+\varepsilon_x\,e^{\,-\,j\,\omega\,h}\,\right
      \rbrace$, i. e.:
\noindent The wavenumber of the finite difference scheme
(\ref{approx}) is thus:

\begin{equation}
 \overline{\lambda}=-\,j\, \left \lbrace\beta_x \,e^{\,0}
      +\delta_x
      \,e^{\,j\,\omega\,h}+\varepsilon_x\,e^{\,-\,j\,\omega\,h}\,\right \rbrace
\end{equation}

\noindent To ensure that the Fourier transform of the finite
difference scheme is a good approximation of the partial derivative
$ (\, \frac{\partial u}{\partial x}\, )_l$ over the range of waves
with wavelength longer than $4\,h$, the a priori unknowns
coefficients $\beta_x$, $\delta_x$, and $\varepsilon_x$ must be
choosen so as to minimize the integrated error:

\begin{equation}\begin{array}{rcl}
 {\mathcal E} &=&\int_{-\frac{\pi}{2}}^{\frac{\pi}{2}} | \lambda \,h- \overline{\lambda}
 \,h|^2\,d(\lambda \,h)\\
 &=& \int_{-\frac{\pi}{2}}^{\frac{\pi}{2}} | \kappa+j\,h\,\left \lbrace\beta_x \,e^{\,0}
      +\delta_x
      \,e^{\,j\,\kappa}+\varepsilon_x\,e^{\,-\,j\,\kappa}\,\right \rbrace\,|^2\,d(\kappa)
      \end{array}
\end{equation}

\noindent The conditions that ${\mathcal E}$ is a minimum are:

\begin{equation}
\frac{\partial  {\mathcal E}}{\partial \beta_x}=\frac{\partial
 {\mathcal E}}{\partial \delta_x}=\frac{\partial  {\mathcal E}}{\partial \varepsilon_x}=0
\end{equation}

\noindent and provide the following system of linear algebraic
equations:

\begin{equation}\left \lbrace
\begin{array}{rcl}
2 \,\pi \,h\, \beta_x +4 \,(h\,\delta_x +h\,\varepsilon_x -1)&=&0\\
 4 \,h\,\beta_x +\pi \, (2\,
\delta_x -1)&=&0\\
4 \,h\,\beta_x +2 \,\pi \, h\,\varepsilon_x &=&0
\end{array} \right.
\end{equation}

\noindent  which enables us to determine the required values of
$\beta_x$, $\delta_x$, and $\varepsilon_x$:

\begin{equation}
\label{Opt_Values1}\left \lbrace
\begin{array}{rcl}
 \beta_x &=& \beta_x^{opt} \,=\, \frac{\pi }{h\,(\pi ^2-8)}\\
\delta_x &=&\delta_x^{opt}\,=\, \frac{1}{2}-\frac{2}{h\,(\pi ^2-8)}\\
\varepsilon_x &=&\varepsilon_x^{opt}\,=\, -\frac{2}{h\,(\pi ^2-8)}\\
\end{array} \right.
\end{equation}

\section{The Sylvester equation}
\label{Sylv}
\subsection{Matricial form of the finite differences problem}

\begin{theorem}
The problem (\ref{scheme}) can be written under the following
matricial form:
\begin{equation}
\label{Eq} {M_1}\,U +U\,M_2+{\cal{L}}(U)=M_0
\end{equation}

\noindent where $M_1$ and $M_2$ are square matrices respectively
$n_x-1$ by $n_x-1$, $n_t$ by $n_t$, given by:
\begin{equation}
\begin{array}{ccc}
{M_1}= \left (
\begin{array}{ccccc}
\beta & \delta & 0 & \ldots &  0 \\
\varepsilon& \beta & \ddots & \ddots &  \vdots \\
0 &  \ddots & \ddots & \ddots &  0\\
\vdots &  \ddots & \ddots & \beta &  \delta\\
0 &  \ldots & 0 & \varepsilon &  \beta\\
\end{array} \right ) &
 & {M_2}= \left (
\begin{array}{ccccc}
0 & \gamma& 0 & \ldots &  0 \\
\alpha & 0 & \ddots & \ddots &  \vdots \\
0 &  \ddots & \ddots & \ddots &  0\\
\vdots &  \ddots & \ddots & \ddots &  \gamma\\
0 &  \ldots & 0 & \alpha &  0\\
\end{array} \right )
\end{array}
  \end{equation}

\noindent the matrix $M_0$ being given by:
\bigskip
  \begin{equation}
  \label{M0}
\scriptsize{{M_0}= \left (
\begin{array}{ccccc}
    -\gamma \,u_1^0
      -\varepsilon \,u_{0}^1-\eta\,u_0^0-\theta\,u_{0}^{2}-\vartheta\,u_{2}^{0}
            &   -\varepsilon \,u_{0}^2 -\eta\,u_0^1-\theta\,u_{0}^{3}& \ldots & \ldots   &
            -\varepsilon \,u_{0}^{n_t}-\eta\,u_{0}^{n_t-1} \\
 -\gamma \,u_2^0-\eta\,u_{1}^{0}-\vartheta\,u_{3}^{0}
       &  0
 & \ldots  & \ldots    &  0 \\
\vdots  &  \vdots & \vdots & \vdots &  \vdots\\
 -\gamma \,u_{n_x-2}^0-\eta\,u_{n_x-2}^{0}-\vartheta\,u_{n_x-1}^{0}
       &  0
 & \ldots  & \ldots   &  0 \\
    -\gamma \,u_{n_x-1}^0-\delta \,u_{n_x}^{1}-\eta\,u_{n_x-2}^{0}-\zeta\,u_{n_x}^{2}-\vartheta\,u_{n_x}^{0}
      &   -\delta \,u_{n_x}^{2} -\zeta\,u_{n_x}^{3}-\vartheta\,u_{n_x}^{1}& \ldots
       & \ldots &  -\delta \,u_{n_x}^{n_t}-\vartheta\,u_{n_x}^{n_t-1}\\
\end{array} \right )}
  \end{equation}

\bigskip
\noindent and where ${\cal {L}}$  is a linear matricial operator
which can be written as: \begin{equation} {\cal {L}}={\cal
{L}}_1+{\cal {L}}_2+{\cal {L}}_3+{\cal {L}}_4
\end{equation}
\noindent where ${\cal {L}}_1$, ${\cal {L}}_2$, ${\cal {L}}_3$ and
${\cal {L}}_4$ are given by:
\begin{equation}
\begin{array}{ccc}
 {\cal {L}}_1(U)  = \zeta  \left (
\begin{array}{ccccc}
u_2^2 & u_2^3 & \ldots & u_2^{n_t} &  0\\
u_3^2 & u_3^3 & \ldots & \vdots &   \vdots \\
\vdots &  \vdots & \ddots & \vdots &  \vdots\\
u_{n_x-1}^2 & u_{n_x-1}^3  & \ldots &  u_{n_x-1}^{n_t}& 0  \\
0 & 0  & \ldots &  0& 0  \\
\end{array}
 \right )
  &
  &  {\cal {L}}_2(U) = \eta \left (
\begin{array}{ccccc}
0 & 0  & \ldots &  0& 0 \\
0 & u_1^1  & u_1^2 & \ldots &  u_1^{n_t-1} \\
0 & u_1^0  & u_1^1 & \ldots &  u_2^{n_t-1} \\
\vdots &  \vdots & \vdots &  \ddots & \vdots \\
0 & u_{n_x-2}^1  & u_{n_x-2}^2 & \ldots &  u_{n_x-2}^{n_t-1} \\
\end{array}
 \right )
\end{array}
  \end{equation}

\begin{equation}
\begin{array}{ccc}
{\cal {L}}_3(U) = \theta \left (
\begin{array}{ccccc}
0 & \ldots &  \ldots &    \ldots &  0 \\
u_1^2 & u_1^3  & \ldots &  u_1^{n_t}& 0 \\
u_2^2 & u_2^3  & \ldots &  u_2^{n_t}& 0 \\
\vdots &  \vdots & \vdots &  \vdots & \vdots \\
u_{n_x-2}^2 & u_{n_x-2}^3  & \ldots &  u_{n_x-2}^{n_t}& 0 \\
\end{array}
 \right )
  &
  &  {\cal {L}}_4(U) = \vartheta \left (
\begin{array}{ccccc}
0 & u_2^1 & u_2^2 & \ldots  &  u_2^{n_t-1} \\
0 & u_3^1 & u_3^2 & \ldots  &  u_3^{n_t-1} \\
\vdots &  \vdots & \ddots  & \ddots &  \vdots\\
0 & u_{n_x-1}^1& \ldots & \ldots  &  u_{n_x-1}^{n_t-1} \\
0 &  0 & \ldots &   \ldots &  0\\
\end{array}
 \right )
 \end{array}
 \end{equation}
\bigskip

\end{theorem}

\bigskip

\begin{proposition}
\noindent The second member matrix $M_0$ bears the initial
conditions, given for the specific value $n=0$, which correspond to
the initialization process when computing loops, and the boundary
conditions, given for the specific values $i=0$, $i=n_x$.
\end{proposition}

\bigskip

\noindent Denote by $u_{exact}$ the exact solution of (\ref{transp}).\\
\noindent The corresponding matrix $U_{exact}$ will be:

\begin{equation}
U_{exact}=[{{{U_{{exact}_i}}}^n}{]_{\, 1\leq i\leq {n_x-1},\, 1\leq
n\leq {n_t}\, }} \end{equation} where:

\begin{equation}
{U_{exact}}_i^n=U_{exact}(x_i,t_n)
 \end{equation}

\noindent with $x_i=i \; h$, $t_n=n \; \tau$. \bigskip

\bigskip

  \begin{definition}
\noindent We will call \textit{error matrix} the matrix defined by:
 \begin{equation}
 \label{err}
E=U-U_{exact}
   \end{equation}
  \end{definition}

\bigskip

\noindent Consider the matrix $F$ defined by:
\begin{equation} F={M_1}\,U_{exact}+U_{exact}\,M_2 + {\cal{L}}(U_{exact})-M_0\end{equation}

\bigskip

  \begin{proposition}
\noindent The \textit{error matrix} $E$ satisfies:

   \begin{equation}
   \label{eqmtr}
{M_1}\,E+E\,M_2+{\cal{L}}(E)=F
   \end{equation}

 \end{proposition}

\subsection{The matrix equation}

   \begin{theorem}
\noindent Minimizing the error due to the approximation induced by
the numerical scheme is equivalent to minimizing the norm of the
matrices $E$ satisfying (\ref{eqmtr}).
   \end{theorem}

   \bigskip

{\em Note:} \noindent Since the linear matricial operator
${\cal{L}}$ appears only in the Crank-Nicolson scheme, we will
restrain our study to the case ${\cal{L}}=0$. The generalization to
the case ${\cal{L}} \neq 0$ can be easily deduced.

   \bigskip

   \begin{proposition}

   \noindent The problem is then the determination of the minimum norm solution
of:

   \begin{equation}
   \label{SylvErr}
{M_1}\,E+E\,M_2=F
   \end{equation}

\noindent which is a specific form of the Sylvester equation:

   \begin{equation}
   \label{SylvGen}
AX+XB=C
   \end{equation}
where $A$ and $B$ are respectively $m$ by $m$ and $n$ by $n$
matrices, $C$ and $X$, $m$ by $n$ matrices.

   \end{proposition}

\bigskip

\subsection{Minimization of the error} \label{MinErr}

\subsubsection{Theory}

\noindent Calculation yields:

\footnotesize \noindent \begin{equation}\left \lbrace
\begin{array}{ccc}
{M_1}\,^T M_1&=& diag \big
 (\left ( \begin{array}{cc}
\beta^2+ \delta^2 & \beta\,(\delta+\varepsilon) \\
\beta\,(\delta+\varepsilon) &\varepsilon^2+ \beta^2 \\
\end{array} \right ),\ldots, \left ( \begin{array}{cc}
\beta^2+ \delta^2 & \beta\,(\delta+\varepsilon) \\
\beta\,(\delta+\varepsilon) &\varepsilon^2+ \beta^2 \\
\end{array} \right )
 \big )\\ {M_2}\,^T M_2&=&
 diag \big
 (\left ( \begin{array}{cc}
\gamma^2 & 0 \\
0 &\alpha^2\\
\end{array} \right ),\ldots, \left ( \begin{array}{cc}
\gamma^2 & 0 \\
0 &\alpha^2 \\
\end{array} \right )
\end{array} \right.
  \end{equation}

\normalsize \noindent The singular values of $M_1$ are the singular
values of the block matrix $\big
 (\left ( \begin{array}{cc}
\beta^2+ \delta^2 & \beta\,(\delta+\varepsilon) \\
\beta\,(\delta+\varepsilon) &\varepsilon^2+ \beta^2 \\
\end{array} \right )$, i. e. \begin{equation} \frac{1}{2} \,(2 \beta ^2+\delta
^2+\varepsilon
   ^2-(\delta +\varepsilon ) \,\sqrt{4 \beta ^2+\delta
   ^2+\varepsilon ^2-2 \delta \, \varepsilon })  \end{equation} \noindent of order $\frac {n_x-1}{2}$, and \begin{equation} \frac{1}{2} \,(2 \beta
^2+\delta ^2+\varepsilon
   ^2+(\delta +\varepsilon ) \,\sqrt{4 \beta ^2+\delta
   ^2+\varepsilon ^2-2 \delta \, \varepsilon })  \end{equation} \noindent of order $\frac {n_x-1}{2}$.\\

\noindent The singular values of $M_2$ are $\alpha^2 $, of order
$\frac {n_t}{2}$, and $\gamma^2 $, of order $\frac {n_t}{2}$.

\noindent Consider the singular value decomposition of the matrices
$M_1$ and $M_2$:

\begin{equation}
U_1^T\,M_1\,V_1=\left (
\begin{array}{cc}
\widetilde{M_1} &0 \\
0& 0\\
\end{array} \right )
\,\,\, ,  \,\,\, U_2^T\,M_1\,V_2=\left (
\begin{array}{cc}
\widetilde{M_2} &0 \\
0& 0\\
\end{array} \right )
  \end{equation}

\noindent where $U_1$, $V_1$, $U_2$, $V_2$, are orthogonal matrices.
\noindent $\widetilde{M_1}$, $\widetilde{M_2}$ are diagonal
matrices, the diagonal terms of which are respectively the nonzero
eigenvalues of the symmetric matrices $M_1\,^T M_1$, $M_2\,^T
M_2$.\\

\noindent Multiplying respectively \ref{SylvErr} on the left side by
$^T U_1$, on the right side by $V_2$, yields:

\begin{equation}
U_1^T\,M_1\,E\,V_2+U_1^T\,E\,M_2\,V_2=U_1^T\,F\,V_2
  \end{equation}

\noindent which can also be taken as:
\begin{equation}
^T U_1\,M_1\,V_1 \,^T V_1\,E\,V_2+^T U_1\,E\,^T U_2\,^T
U_2\,M_2\,V_2=U_1^T\,F\,V_2
  \end{equation}

\noindent Set:

\begin{equation}
^T V_1\,E\,V_2= \left (
\begin{array}{cc}
\widetilde{E_{11}} & \widetilde{E_{12}} \\
\widetilde{E_{21}} & \widetilde{E_{22}}\\
\end{array} \right )
\,,\,^T U_1\,E\,^T U_2= \left (\begin{array}{cc}
\widetilde{\widetilde{E_{11}}} & \widetilde{\widetilde{E_{12}}} \\
\widetilde{\widetilde{E_{21}}} & \widetilde{\widetilde{E_{22}}}\\
\end{array} \right )
  \end{equation}

\begin{equation}
\label{Ftilde} ^T U_1\,F\,V_2= \left (
\begin{array}{cc}
\widetilde{F_{11}} & \widetilde{F_{12}} \\
\widetilde{F_{21}} & \widetilde{F_{22}}\\
\end{array} \right )
  \end{equation}
\noindent We have thus:
\begin{equation}\left (\begin{array}{cc}
\widetilde{M_{1}}\,\widetilde{E_{11}} & \widetilde{M_{1}}\,\widetilde{E_{12}} \\
0 & 0\\
\end{array} \right )+\left (\begin{array}{cc}
\widetilde{\widetilde{E_{11}}}\,\widetilde{M_{2}} & 0 \\
\widetilde{\widetilde{E_{21}}}\,\widetilde{M_{2}} & 0\\
\end{array} \right )=\left (\begin{array}{cc}
\widetilde{F_{11}} & \widetilde{F_{12}} \\
\widetilde{F_{21}} & \widetilde{F_{22}}\\
\end{array} \right )
  \end{equation}

\noindent It yields:
\begin{equation}
\left \lbrace
\begin{array}{ccc}
\widetilde{M_{1}}\,\widetilde{E_{11}} +
\widetilde{\widetilde{E_{11}}}\,\widetilde{M_{2}}&=&\widetilde{F_{11}}\\
\widetilde{M_{1}} \,\widetilde{E_{12}}&=&\widetilde{F_{12}}\\
\widetilde{\widetilde{E_{21}}}\,\widetilde{M_{2}}&=&\widetilde{F_{21}}\\
\end{array} \right.
  \end{equation}

\noindent One easily deduces:
\begin{equation}
\left \lbrace
\begin{array}{ccc}
\widetilde{E_{12}}&=&{\widetilde{M}_{1}}^{-1}\,\widetilde{F_{12}}\\
{\widetilde{\widetilde{E}_{21}}}&=&\widetilde{F_{21}}\,{\widetilde{M_{2}}}^{-1}\\
\end{array} \right.
  \end{equation}

\noindent The problem is then the determination of the
$\widetilde{E_{11}}$ and $\widetilde{\widetilde{E_{11}}}$
satisfying:

\begin{equation}
\label{Pb} \widetilde{M_{1}}\,\widetilde{E_{11}} +
\widetilde{\widetilde{E_{11}}}\,\widetilde{M_{2}}=\widetilde{F_{11}}
  \end{equation}

\noindent Denote respectively by $\widetilde{e_{ij}}$,
$\widetilde{\widetilde{e_{ij}}}$ the components of the matrices
$\widetilde{E}$, $\widetilde{\widetilde{E}}$.\\
\noindent The problem \ref{Pb} uncouples into the independent
problems:\\ \noindent minimize
\begin{equation}
\sum_{i,j} {\widetilde{e_{ij}}}^2+{\widetilde{\widetilde{e_{ij}}}}^2
\end{equation}

 \noindent under the constraint \begin{equation}
\widetilde{M_{1}}_{ii}\,
{\widetilde{e_{ij}}}+\widetilde{M_{2_{ii}}}\,{\widetilde{\widetilde{e_{ij}}}}=\widetilde{F_{11}}_{ij}
  \end{equation}
\noindent This latter problem has the solution:

\begin{equation}\left \lbrace
\begin{array}{ccc}
\widetilde{e_{ij}}
&=&\frac{\widetilde{{M_{1}}_{ii}}\,\widetilde{{F_{11}}_{ij}}}
{{\widetilde{{M_{1}}_{ii}}}^2+{\widetilde{{M_{2}}_{jj}}^2}}\\
\widetilde{\widetilde{e_{ij}}} &=&\frac{\widetilde{{M_{2}}
_{jj}}\,\widetilde{{F_{11}}_{ij}}}{{\widetilde{{M_{1}}_{ii}}}^2+{\widetilde{{M_{2}}_{jj}}^2}}\\
\end{array} \right.
  \end{equation}

\noindent The minimum norm solution of \ref{SylvErr} will then be
obtained when the norm of the matrix $\widetilde{{F_{11}}}$ is
minimum.\\
\noindent In the following, the euclidean norm will be considered.

\noindent Due to (\ref{Ftilde}):
 \begin{equation}
  \|\widetilde{{F_{11}}}  \| \leq  \|\widetilde{{F}}  \| \leq \|U_1
  \| \,\|F \|\,\|V_2  \| \leq \|U_1
  \| \,\|V_2  \| \, \|M_1\,U_{exact}+U_{exact}\,M_2-M_0
  \|
   \end{equation}

\noindent $U_1$ and $V_2$ being orthogonal matrices, respectively
$n_x-1$ by $n_x-1$, $n_t$ by $n_t$, we have:

 \begin{equation}
  \|U_1
  \|^2 =n_x-1\,\,\,, \,\,\,\|V_2
  \|^2 =n_t
   \end{equation}

\noindent Also:

 \begin{equation}
  \|M_1
  \|^2 =\frac{n_x-1}{2}\,\big ( 2\,\beta^2+\delta^2+\varepsilon^2 \big )\,\,\,, \,\,\,
  \|M_2
  \|^2 =\frac{n_t}{2}\,\big ( \alpha^2+\gamma^2 \big )
   \end{equation}

\noindent The norm of $M_0$ is obtained thanks to relation
(\ref{M0}).

\noindent This results in:
 \begin{equation}
 \label{Min}
  \|\widetilde{{F_{11}}} \| \leq \sqrt {n_t\,(n_x-1)} \, \left \lbrace  \| U_{exact} \| \,\big (
  \sqrt{\frac{n_x-1}{2}}\,\sqrt{ 2\,\beta^2+\delta^2+\varepsilon^2 }+
   \sqrt{\frac{n_t}{2}}\,\sqrt{\alpha^2+\gamma^2 } \,\big )+
  \|M_0  \|\right \rbrace
   \end{equation}

\noindent $ \|\widetilde{{F_{11}}} \| $ can be minimized through the
minimization of the second factor of the right-side member of
(\ref{Min}), which is function of the scheme parameters.\\

\noindent $\| U_{exact} \| $ is a constant. The quantities $
\sqrt{\frac{n_x-1}{2}}\,\sqrt{ 2\,\beta^2+\delta^2+\varepsilon^2
  }$, $\sqrt{\alpha^2+\gamma^2 } $ and $\|M_0  \|$ being
  strictly positive, minimizing the second factor of the right-side member
  of (\ref{Min}) can be obtained through the minimization of the
  following functions:

 \begin{equation}
 \label{func_Sylv}
 \left \lbrace
\begin{array}{rcl}
f_1(\beta,\delta,\varepsilon)&=&\sqrt{ 2\,\beta^2+\delta^2+\varepsilon^2}\\
 f_2(\alpha,\gamma)&=&\sqrt{\alpha^2+\gamma^2 }\\
 f_3(\alpha,\beta,\gamma,\delta,\varepsilon)&=&\|M_0  \|\\
\end{array} \right.
\end{equation}

\section{Numerical example: a new DRP scheme}

\label{Ex}

\noindent Consider the scheme (\ref{scheme}) where the values of
$\beta_x$, $\delta_x$, and $\varepsilon_x$ are given by
(\ref{Opt_Values1}).\\

\noindent Let, in a first time, the values of the coefficients
$\alpha$, $\beta_t$, $\gamma$, $\delta_t$, and $\varepsilon_t$
remain unknown, and advect a sinusoidal signal
 \begin{equation}  u=\cos\,[\,k \,(x-c\,t)\,]
\end{equation}
\noindent through this scheme, with Dirichlet boundary
conditions. ($c$ is taken equal to 1, and $k= {\pi}$).\\

\noindent Calculation yields then:

 \begin{equation}\scriptsize{
 \left \lbrace
\begin{array}{rcl}
f_1(\beta,\delta_x^{opt}+\delta_t,\varepsilon^{opt}+\varepsilon_t)&=&\sqrt{2
\left(\,{\beta_t}+\frac{\pi }{4 \,h
   \left(\pi
   ^2-8\right)}\right)^2+\left(\,{\delta_t}+\frac{\frac{1}{2}-\frac{2}{\pi
   ^2-8}}{2
   h}\right)^2+\left(\,{\varepsilon_t}+\frac{\frac{1}{2}-\frac{2}{\pi
   ^2-8}}{2
   h}\right)^2}\\
 f_2(\alpha,\gamma)&=&\sqrt{\alpha^2+\gamma^2 }\\
 f_3(\alpha,\beta_x^{opt}+\beta_t,\gamma,\delta_x^{opt}+\delta_t,\varepsilon^{opt}+\varepsilon_t)&=&\sqrt{3 \gamma ^2+3
\left(\,{\delta_t}+\frac{\frac{1}{2}-\frac{2}{-8+\pi ^2}}{2
   h}\right)^2+\left(\gamma -\,{\varepsilon_t}-\frac{\frac{1}{2}-\frac{2}{\pi
   ^2-8}}{2
   h}\right)^2+3 \left(\,{\varepsilon_t}+\frac{\frac{1}{2}-\frac{2}{\pi
   ^2-8}}{2
   h}\right)^2}\\
\end{array} \right.}
\end{equation}

\noindent Minimum values for $f_1$ and $f_3$ can thus be obtained
choosing negative values for $\beta_t$, while choosing positive ones
for $\delta_t$ and $\varepsilon_t$, the absolute values of which are
respectively close to those of $\beta_x$, $\delta_x$ and
$\varepsilon_x$. $f_2$ is minimized choosing $\gamma=0$.\\

\noindent In the following, we have choosen to set:

 \begin{equation}
 \left \lbrace
\begin{array}{rcl}
\beta_t&=&-0.9\,\beta_x^{opt}\\
\delta_t&=&-0.9\,\delta_x^{opt}\\
\varepsilon_t&=&-0.9\,\varepsilon_x^{opt}\\
\end{array} \right.
\end{equation}

\noindent and $\alpha=10$.

\begin{figure}[ht]
\hspace{4cm}
\begin{tabular}{c}
\psfig{height=4.5cm,width=7cm,angle=0,file=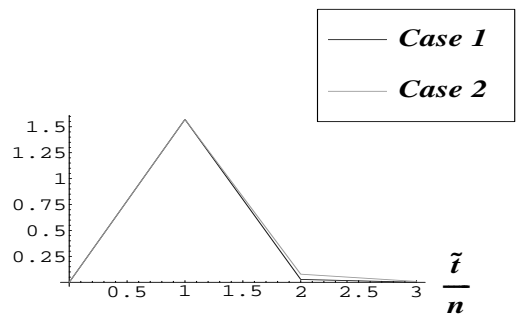}\\
\end{tabular}
\caption{\small{Value of the $L_2$ norm  of the error}. }
\label{Opt}
\end{figure}

\noindent The value of the $L_2$ norm of the error, for:
\begin{enumerate}
\item [\emph{i}.] case 1: our new scheme, with $cfl=0.9$;
\item [\emph{ii}.] case 2: the Lax scheme, with $cfl=0.9$;
\end{enumerate}
\noindent is displayed in Figure \ref{Opt}. The error curve
corresponding to the first case is the minimal one.

\newpage

\section{Conclusion}

The above results open new ways for the building of DRP schemes. It
seems that the research on this problem has not been performed
before as far as our knowledge goes. In the near future, we are
going to extend the techniques described herein to nonlinear
schemes, in conjunction with other innovative methods as the Lie
group theory.

\addcontentsline{toc}{section}{\numberline{}References}

\end{document}